\documentclass[a4paper,10pt]{amsart}
\usepackage{amsmath,amsthm,latexsym,amscd,amssymb}
\usepackage[all]{xy}
\usepackage{graphics}
\makeindex
\numberwithin{equation}{subsection}

\newtheorem{prop}{Proposition}[section]
\newtheorem{lem}[prop]{Lemma}
\newtheorem{ddd}[prop]{Definition}
\newtheorem{theorem}[prop]{Theorem}
\newtheorem{cor}[prop]{Corollary}

\newcommand{\ind}{\mathop{\mbox{\rm ind}}}
\newcommand{\dom}{\mathop{\rm dom}}

\newcommand{\Tr}{\mathop{\rm Tr}}
\newcommand{\U}{{\mathcal U}}
\newcommand{\C}{C^{\infty}}
\newcommand{\ve}{\varepsilon}
\DeclareMathOperator{\supp}{supp}
\DeclareMathOperator{\spfl}{sf}

\DeclareMathOperator{\Ker}{Ker}

\def\bbbr{{\rm I\!R}} 
\def\bbbn{{\rm I\!N}} 

\def\bbbc{{\rm I\!C}}

\def\bbbq{{\mathchoice {\setbox0=\hbox{$\displaystyle\rm Q$}\hbox{\raise
0.15\ht0\hbox to0pt{\kern0.4\wd0\vrule height0.8\ht0\hss}\box0}}
{\setbox0=\hbox{$\textstyle\rm Q$}\hbox{\raise
0.15\ht0\hbox to0pt{\kern0.4\wd0\vrule height0.8\ht0\hss}\box0}}
{\setbox0=\hbox{$\scriptstyle\rm Q$}\hbox{\raise
0.15\ht0\hbox to0pt{\kern0.4\wd0\vrule height0.7\ht0\hss}\box0}}
{\setbox0=\hbox{$\scriptscriptstyle\rm Q$}\hbox{\raise
0.15\ht0\hbox to0pt{\kern0.4\wd0\vrule height0.7\ht0\hss}\box0}}}}

\def\bbbz{{\mathchoice {\hbox{$\sf\textstyle Z\kern-0.4em Z$}}
{\hbox{$\sf\textstyle Z\kern-0.4em Z$}}
{\hbox{$\sf\scriptstyle Z\kern-0.3em Z$}}
{\hbox{$\sf\scriptscriptstyle Z\kern-0.2em Z$}}}}

\def\bbbc{{\mathchoice {\setbox0=\hbox{$\displaystyle\rm C$}\hbox{\hbox
to0pt{\kern0.4\wd0\vrule height0.9\ht0\hss}\box0}}
{\setbox0=\hbox{$\textstyle\rm C$}\hbox{\hbox
to0pt{\kern0.4\wd0\vrule height0.9\ht0\hss}\box0}}
{\setbox0=\hbox{$\scriptstyle\rm C$}\hbox{\hbox
to0pt{\kern0.4\wd0\vrule height0.9\ht0\hss}\box0}}
{\setbox0=\hbox{$\scriptscriptstyle\rm C$}\hbox{\hbox
to0pt{\kern0.4\wd0\vrule height0.9\ht0\hss}\box0}}}}

\begin{document}

\title{Spectral flow as winding number and integral formulas}

\author{Charlotte Wahl}

\email{ac.wahl@web.de}



\begin{abstract}
A general integral formula for the spectral flow of a path of unbounded selfadjoint Fredholm operators subject to certain summability conditions is derived from the interpretation of the spectral flow as a winding number.
\end{abstract}

\maketitle

\section*{Introduction}

Integral formulas for the spectral flow play an important role in recent developments in noncommutative geometry (see \cite{bcp} and references therein). Here we have in mind formulas of the form 
$$\spfl((D_t)_{t \in [0,1]})=\int_0^1(\frac{d}{dt}D_t) \psi(D_t)~dt$$  for a path $(D_t)_{t \in [0,1]}$ of unbounded selfadjoint Fredholm operators such that certain summability and differentiability conditions are fulfilled. For the moment we also assume that $D_0$ and $D_1$ are invertible and unitarily equivalent. If this is not the case, there is an additional contribution from the endpoints. 

Integral formulas have been proven for bounded perturbations and special classes of functions $\psi$, most importantly for (up to normalization) $\psi(x)=(1+x^2)^{-p/2}$ with $p>0$ and $\psi(x)=e^{-tx^2}$ for $t>0$. These have applications to $p$-summable and $\theta$-summable Fredholm modules respectively. 
Getzler suggested and proved the formula in the $\theta$-summable case \cite{gl}, and Carey--Phillips in the $p$-summable case \cite{cp1}. In both cases the formulas have been generalized and extended to Breuer--Fredholm operators by Carey--Phillips \cite{cp1}\cite{cp2}. 

In this paper we prove a general integral formula using the interpretation of the spectral flow as a winding number. In contrast to loc.cit. our proof works also for $D_t-D_0$ unbounded. 
As applications we discuss integral formulas for paths of elliptic operators on a closed manifold and review the $p$-summable case and the $\theta$-summable case.

{\it Acknowledgement:} I thank the anonymous referee of a previous version of this paper for useful suggestions.

\section{Derivation of the integral formula}

The spectral flow of a norm continuous path of bounded selfadjoint Fredholm operators was introduced in \cite{aps}. It measures the net number of eigenvalues changing sign along the path. There are several suitable topologies on the space of unbounded selfadjoint Fredholm operators for which the spectral flow is well-defined. See \cite{l} for a comparison. Recall that a path is gap continuous if the resolvents depend continuously on the parameter. In \cite{blp} the spectral flow of a gap continuous path was defined and expressed in terms of the winding number of its Cayley transform. The definition of the spectral flow in terms of the winding number that we will give below is closely related and is justified by \cite[Prop. 2.6]{wa}.  

Let $H,H'$ be separable Hilbert spaces. Let $B(H,H')$ be the space of bounded operators from $H$ to $H'$ endowed with the norm topology. Let $K(H) \subset B(H)$ be the ideal of compact operators. The group of unitaries in $B(H)$ is denoted by ${\mathcal U}(H)$.

Let $${\mathcal U}_K(H)=\{U \in {\mathcal U}(H) ~|~U-1 \in K(H)\} \ .$$

We denote by $l^p(H)$ for $p \in [1,\infty)$ the Banach ideal of $p$-summable operators in $B(H)$ with norm 
$$\|A\|_p= (\Tr |A|^p)^{\frac 1p}$$
and by $l^{\infty}(H)$ the space of bounded operators on $H$ endowed with the strong operator topology. Furthermore for $p \in [1,\infty]$ let $l_{sa}^p(H) \subset l^p(H)$ be the subspace of selfadjoint operators with the subspace topology. 

Per definition a selfadjoint operator $D$ on $H$ is Fredholm if its bounded transform $D(1+D^2)^{-\frac 12}$ is Fredholm. Every selfadjoint operator with compact resolvents is Fredholm. 

We will use that for a selfadjoint operator $D$ with compact resolvents and a continuous map $[0,1] \to l_{sa}^{\infty}(H),~ t \mapsto A_t$ the maps $$[0,1] \to K(H),~ t \mapsto (D+ A_t \pm i)^{-1}$$ are continuous.

We call a function $\chi\in C^1(\bbbr)$ fulfilling
\begin{itemize}
\item $\chi^{-1}(0)=\{0\}$,
\item $\lim_{x \to -\infty}\chi(x)=-1$ and $\lim_{x \to  \infty}\chi(x)=1$, 
\item $\chi' \in C_0(\bbbr)$ with $\chi'(0)>0$ and $\chi' \ge 0$
\end{itemize}  
a normalizing function. (This definition differs slightly from the one in \cite{wa}). 

If $(D_t)_{t \in [0,1]}$ is a gap continuous path of selfadjoint operators with compact resolvents, then $f(D_t)$ is compact for any $f\in C_0(\bbbr)$ and continuous in $t$. In particular, if $\chi$ is a normalizing function, then $\chi(D_t)^2-1$ is compact and depends continuously on $t$; furthermore since $e^{\pi i(\chi+1)}-1 \in C_0(\bbbr)$, the unitary $e^{\pi i(\chi(D_t)+1)} \in {\mathcal U}_K(H)$ depends continuously on $t$ as well.

Let $S^1=[0,1]/_{0\sim 1}$.

The winding number of a loop $s:S^1 \to {\mathcal U}_K(H)$ with $s-1 \in C^1(S^1,l^1(H))$ is defined by $$w(s)= \frac{1}{2\pi i}\int_0^1 \Tr(s(x)^{-1} s'(x))~dx \ .$$ The formula holds also for loops that are piecewise in $C^1$ (see \cite{kl}). The winding number extends to a well-defined isomorphism $$w:\pi_1({\mathcal U}_K(H)) \cong \bbbz \ .$$

We will work with the following definition of the spectral flow. We only need the case of operators with compact resolvents. See \cite{wa} for a general definition.

\begin{ddd}
Let $(D_t)_{t \in [0,1]}$ be a gap continuous path of selfadjoint operators with compact resolvents. Assume that $D_0$ and $D_1$ are invertible. Let $\chi$ be a normalizing function such that $\chi(D_0)$ and $\chi(D_1)$ are involutions.  Then the spectral flow of $(D_t)_{t \in [0,1]}$ is given by
$$\spfl((D_t)_{t \in [0,1]})=w([e^{\pi i(\chi(D_t)+1)}]) \ .$$
\end{ddd}

The definition does not depend on the choice of the normalizing function. Furthermore the spectral flow is homotopy invariant, invariant under orientation preserving reparametrisation and additive with respect to concatenation of paths. It can be defined for arbitrary intervals via reparametrisation. Sometimes we write $\spfl(D_t)$ for $\spfl((D_t)_{t \in [0,1]})$.

The definition can be extended to the case where $D_0$, $D_1$ are not necessarily invertible: Choose a real-valued function $\phi \in \C([-1,2])$ supported in $[-1,-\frac 12]\cup [\frac 32,2]$, with $\phi(-1)=\phi(2)$ and such that $[-\phi(-1),0)$ is a subset of the intersection of the resolvent sets of $D_0$ and $D_1$.  The spectral flow of $(D_t)_{t \in [0,1]}$ is defined as the spectral flow of the path $(\tilde D_t)_{t \in [-1,2]}$ given by $\tilde D_t:=D_0+\phi(t)$ for $t \in [-1,0]$,~$\tilde D_t=D_t$ for $t \in [0,1]$ and $\tilde D_t=D_1+\phi(t)$ for $t \in [1,2]$. The definition is independent of the choices. 

If $(D_t)_{t \in [0,1]}$ is a loop, then $$\spfl((D_t)_{t \in [0,1]})=w([e^{\pi i(\chi(D_t)+1)}])$$
for any normalizing function $\chi$.

For simplicity we will only consider paths with invertible endpoints. Using the definition one can then easily deduce formulas for the general case. 

In the following lemma, in preparation for the general case, we derive an integral formula for the spectral flow on a finite-dimensional Hilbert space. In its proof we will use some facts about the relative index of projections (see \cite{ass}): 

The relative index of a pair $(P,Q)$ of projections on $H$ with $P-Q \in K(H)$ is defined by $$\ind(P,Q)=\ind(QP:P(H) \to Q(H)) \ .$$ 
If $P-Q \in l^1(H)$, then $$\ind(P,Q)=\Tr(P-Q) \ .$$
Furthermore $$\spfl(t(2P-1)+(1-t)(2Q-1))=\ind(P,Q) \ .$$

\begin{lem}
\label{findim}

Assume that $H=\bbbc^n$.
Let $(D_t)_{t \in [0,1]}$ be a path of selfadjoint matrices with $(t \mapsto D_t) \in C^1([0,1],M_n(\bbbc))$ and assume that $D_0,D_1$ are invertible.  

Let $\chi \in C^1(\bbbr)$ be a normalizing function. Then 
\begin{eqnarray*}
\spfl((D_t)_{t \in [0,1]})&=&\frac 12 \int_0^1  \Tr (\frac{d}{dt}D_t)\chi'(D_t)~ dt\\
&& +~ \frac 12 \Tr(2P_1-1-\chi(D_1)) - \frac 12 \Tr(2P_0-1-\chi(D_0))\ ,
\end{eqnarray*} 
where $P_i= 1_{\ge 0}(D_i)$.
\end{lem}

\begin{proof}
For the moment assume that the restriction of $\chi$ to $[-c,c]$ is a polynomial, where $c>0$ is such that the spectrum of $D_t$ is contained in $[-c,c]$ for all $t \in [0,1]$.

Set $F_t=\chi(D_1) + t(\chi(D_0)-\chi(D_1))$ and denote
$$L=\frac{1}{2 \pi i} \int_0^1 \Tr  e^{-\pi i(F_t+1)}\frac{d}{dt}e^{\pi i (F_t+1)}~dt - \spfl((F_t)_{t \in [0,1]}) \ .$$
Define a loop $(Q_t)_{t \in [0,2]}$ by $Q_t:= \chi(D_t)$ for $t \in [0,1]$ and $Q_t:=F_{t-1}$ for $t \in [1,2]$. 

Using the definition of the spectral flow we get  
\begin{eqnarray*}
\spfl((Q_t)_{t \in [0,2]}) &=& \frac{1}{2 \pi i} \int_0^1 \Tr  e^{-\pi i(\chi(D_t)+1)}\frac{d}{dt}e^{\pi i ( \chi(D_t)+1)}~dt\\
&&+~ \frac {1}{2 \pi i} \int_0^1 \Tr  e^{-\pi i (F_t+1)}\frac{d}{dt}e^{\pi i (F_t+1)}~dt \ .
\end{eqnarray*}
Since $\spfl((D_t)_{t \in [0,1]})+ \spfl((F_t)_{t \in [0,1]})=\spfl((Q_t)_{t \in [0,2]})$, it follows that
$$\spfl((D_t)_{t \in [0,1]})= \frac{1}{2 \pi i} \int_0^1 \Tr  e^{-\pi i(\chi(D_t)+1)}\frac{d}{dt}e^{\pi i( \chi(D_t)+1)}~dt + L \ . $$
The function $g:[-c,c]\to \bbbc,~g(\lambda)=e^{\pi i(\chi(\lambda)+1)}$ extends to an entire function. Hence if $\Gamma$ is a closed curve in $\bbbc$ not intersecting $[-c,c]$ and with winding number $1$ with respect to the origin, we have that
$$g(D_t)= \frac{1}{2 \pi i}\int_{\Gamma}g(\lambda) (D_t-\lambda)^{-1} d\lambda \ .$$
From $$(D_{t+h}-\lambda)^{-1}-(D_t-\lambda)^{-1}=(D_{t+h}-\lambda)^{-1}(D_t-D_{t+h})(D_t-\lambda)^{-1}$$ we conclude that $$\frac{d}{dt}(D_t-\lambda)^{-1}=-(D_t-\lambda)^{-1}(\frac{d}{dt}D_t)(D_t-\lambda)^{-1} \ ,$$
hence
$$\frac{d}{dt}g(D_t)= -\frac{1}{2 \pi i}\int_{\Gamma}g(\lambda) (D_t-\lambda)^{-1}(\frac{d}{dt}D_t)(D_t-\lambda)^{-1} d\lambda \ .$$
It follows that
\begin{eqnarray*}
\lefteqn{\frac{1}{2 \pi i}\int_0^1 \Tr  e^{-\pi i(\chi(D_t)+1)}\frac{d}{dt}e^{\pi i(\chi(D_t)+1)}~dt}\\
&=& -  \left(\frac{1}{2 \pi i}\right)^2\int_0^1 \Tr  \int_{\Gamma}g(\lambda) e^{-\pi i(\chi(D_t)+1)} (D_t-\lambda)^{-1}(\frac{d}{dt}D_t)(D_t-\lambda)^{-1} d\lambda ~dt \\
&=&  - \left(\frac{1}{2 \pi i}\right)^2 \int_0^1  \Tr \bigl((\frac{d}{dt}D_t)e^{-\pi i(\chi(D_t)+1)} \int_{\Gamma}g(\lambda)  (D_t-\lambda)^{-2} d\lambda\bigr) ~dt\\
&=& \frac{1}{2 \pi i} \int_0^1  \Tr (\frac{d}{dt}D_t)e^{-\pi i(\chi(D_t)+1)}g'(D_t) ~dt\\
&=& \frac{1}{2}  \int_0^1  \Tr (\frac{d}{dt}D_t)\chi'(D_t) ~dt \ .
\end{eqnarray*}

We evaluate $L$: An analogous calculation shows that
$$\frac{1}{2\pi i}\int_0^1 \Tr  e^{-\pi i (F_t+1)}\frac{d}{dt}e^{\pi i (F_t+1)}~dt= \frac{1}{2}\Tr(F_1-F_0)=\frac{1}{2}\Tr(\chi(D_0)-\chi(D_1)) \ .$$

The path $F_t=\chi(D_1) + t(\chi(D_0)-\chi(D_1))$ is homotopic to the path $(1-t)(2P_1-1)+t(2P_0+1)$ by a homotopy of paths with invertible endpoints. Hence, by homotopy invariance of the spectral flow and the properties of the relative index of projections, 
\begin{eqnarray*}
\spfl((F_t)_{t \in [0,1]})&=&\spfl((1-t)(2P_1-1)+t(2P_0+1))\\
&=&\Tr(P_0-P_1)\\
&=&\frac 12\Tr((2P_0-1)-(2P_1-1)) \ .
\end{eqnarray*}
It follows that
\begin{eqnarray*}
L&=&\frac 12\Tr(\chi(D_0)-\chi(D_1))- \frac 12\Tr((2P_0-1)-(2P_1-1)) \\
&=& \frac 12 \Tr((2P_1-1)- \chi(D_1))- \frac 12 \Tr((2P_0-1)- \chi(D_0)) \ .
\end{eqnarray*}

Both sides of the formula are continuous in $\chi$ with respect to the norm of $C^1([-c,c])$. Since we can approximate the restriction of any normalizing function to $[-c,c]$ in $C^1([-c,c])$ by  polynomials, the formula holds for any normalizing function.   
\end{proof}

It is a well-known fact from functional analysis that if $(D_t)_{t \in [0,1]}$ is a gap continuous path of selfadjoint operators, then $[0,1] \to l^{\infty}(H),~ t \mapsto f(D_t)$ is continuous for any $f \in C(\bbbr)$.
Hence if in addition $g \in C(\bbbr)$ such that $[0,1] \to l^1(H),~t \mapsto g(D_t)$ is continuous, then $[0,1] \to l^1(H),~t \mapsto f(D_t)g(D_t)$ is continuous. We will tacitly make use of this property.

For a selfadjoint operator $D$ let $H(D)$ be the Hilbert space whose underlying vector space is $\dom D$ and whose scalar product is given by $<v,w>_D=<v,w>+<Dv,Dw>$, where $<~,~>$ is the scalar product on $H$.

\begin{theorem}
Let $(D_t)_{t \in [0,1]}$ be a path of selfadjoint operators with compact resolvents and common domain such that $(t \mapsto D_t) \in C^1([0,1],B(H(D_0),H))$ and let $D_0,D_1$ be invertible.

Let $\phi \in C_0(\bbbr)$ be such that $\supp \phi$ is connected and $\phi$ is positive on the interior $(\supp \phi)^o$ of its support and such that $[0,1] \to l^1(H),~t \mapsto \phi(D_t)$ is continuous.

Let $\chi \in C^1(\bbbr)$ be a normalizing function with $\supp \chi' \subset (\supp \phi)^o$ and such that there is $C>0$ with $$|\chi^2(x)-1| \le C\phi(x) \mbox{ and }(|x|+1)^p|\chi'(x)| \le C\phi(x) \ ,$$ where $p=0$ if $D_0-D_t$ is bounded for all $t \in [0,1]$ and else $p=1$. 

Then, with $P_i:= 1_{\ge 0}(D_i)$, 
\begin{eqnarray*}
\spfl((D_t)_{t \in [0,1]})&=&\frac 12 \int_0^1  \Tr (\frac{d}{dt}D_t)\chi'(D_t)~ dt\\
&& + ~ \frac 12 \Tr(2P_1-1-\chi(D_1)) - \frac 12 \Tr(2P_0-1-\chi(D_0))\ .
\end{eqnarray*} 

\end{theorem}

Note that each term of the formula is invariant under conjugation by a path of unitaries $(U_t)_{t \in [0,1]}$ such that $(t \mapsto U_t),~(t \mapsto U_t^*) \in C^1([0,1],l^{\infty}(H)) \cap C^1([0,1],l^{\infty}(H(D_0)))$: This is clear for the left hand side and for the contribution of the endpoints. Furthermore   
\begin{eqnarray*}
\lefteqn{\Tr (\frac{d}{dt}U_tD_tU_t^*)\chi'(U_tD_tU_t^*)}\\
&=& \Tr (\frac{d}{dt}U_t)D_tU_t^*\chi'(U_tD_tU_t^*)  + \Tr U_tD_t(\frac{d}{dt}U_t^*)\chi'(U_tD_tU_t^*) \\
&& ~ +~ \Tr U_t(\frac{d}{dt}D_t)U_t^*\chi'(U_tD_tU_t^*) \\
&=&\Tr (\frac{d}{dt}D_t)\chi'(D_t) \ .
\end{eqnarray*}
Here we used the cyclicity of the trace and that $$(\frac{d}{dt}U_t)U_t^* +  U_t(\frac{d}{dt}U_t^*)=\frac{d}{dt}(U_tU_t^*)=0 \ .$$

\begin{proof}
We note that by the remark preceding the theorem the maps  $[0,1] \to l^1(H),~t \mapsto (\chi(D_t)^2-1)$ and $[0,1] \to l^1(H),~t \mapsto (\frac{d}{dt} D_t)\chi'(D_t)$  are well-defined and continuous.

First assume that $\chi' \in \C_c(\bbbr)$ and let $R>0$ be such that $\supp \chi' \in [-R,R]$. We also assume that there are $\lambda,\mu \in C^1([0,1])$ with $\lambda,\mu>R$  such that $\{-\lambda(t),\mu(t)\}$ is a subset of the resolvent set of $D_t$ for all $t \in [0,1]$. Then $P_t=1_{[-\lambda(t),\mu(t)]}(D_t)$ is a path of projections with finite-dimensional range. We have that $D_t=P_tD_t \oplus (1-P_t)D_t$ and $\chi'(D_t)=\chi'(P_tD_t)$. For any $t_0 \in [0,1]$ there is a closed curve $\Gamma$ and $\ve>0$ such that $$P_t=\frac{1}{2\pi i}\int_{\Gamma}(D_t- \lambda)^{-1}d\lambda$$ for all $t \in (t_0-\ve,t_0+\ve)$. Since the resolvents depend differentiably on $t$ in $B(H,H(D_0))$, it follows that $(t \mapsto P_t) \in C^1([0,1],B(H,H(D_0)))$. We may find a family of isometries $(U_t:P_t(H) \to P_0(H))_{t \in [0,1]}$ such that $(t \mapsto U_tP_t) \in C^1([0,1],K(H))$ and such that $U_tP_tU_t^*=P_0$ for each $t \in [0,1]$. 

We define the operator $\tilde D_t:=U_tP_t D_t P_t U_t^*$ on $P_0(H)$. 
Clearly $\spfl((D_t)_{t \in [0,1]})=\spfl((\tilde D_t)_{t \in [0,1]})$.

Lemma  \ref{findim}  applies to $(\tilde D_t)_{t\in[0,1]}$. Furthermore 
\begin{eqnarray*}
\lefteqn{\Tr (\frac{d}{dt}\tilde D_t)\chi'(\tilde D_t)}\\
&=& \Tr \frac{d}{dt}(U_tP_t)D_t\chi'(D_t) P_tU_t^*  + \Tr D_t\frac{d}{dt}(P_t U_t^*)U_tP_t\chi'(D_t)  + \Tr (\frac{d}{dt}D_t)\chi'(D_t) \\
&=&\Tr (\frac{d}{dt}P_t)D_t\chi'(D_t)  + \Tr (\frac{d}{dt}D_t)\chi'(D_t)\\
&=&\Tr (\frac{d}{dt}D_t)\chi'(D_t) \ .
\end{eqnarray*}
Here we used that  
$$\Tr (\frac{d}{dt}P_t)D_t\chi'(D_t)=\Tr P_t(\frac{d}{dt}P_t)P_t D_t\chi'(D_t) =0$$ 
since $P_t(\frac{d}{dt}P_t)P_t=0$. 

This proves the assertion under the two additional assumptions of the beginning of the proof.

We still let $\chi$ and $R$ be as above, but now we drop the assumption on the existence of $\lambda(t),\mu(t)$. We may find a partition $0=x_0< x_1 < \dots < x_{k+1}=1$ such that $\lambda(t),\mu(t)$ as above exist on  $[x_i,x_{i+1}]$ for each $i$. We cannot apply the above arguments directly since $D_{x_i}$ need not be invertible for $i=1, \dots,k$. Let $Q$ be the projection onto $\Ker D_{x_1} + \Ker D_{x_2} + \dots \Ker D_{x_k}$. Let $\ve>0$ be such that $[-\ve,0)$ is in the resolvent set of $D_{x_i}$ for any $i=1, \dots, k$ and let $\phi \in \C([0,1])$ be a real-valued function such that $\phi(1)=\ve$ and $\supp \phi \subset (0, 1]$. For $i=1, \dots, k-1$ we define $(\tilde D^i_t)_{t \in [x_i-1,x_{i+1}+1]}$ by $\tilde D^i_t:=D_{x_i}+\phi(x_i-t)Q$ for $t \in [x_i-1,x_i]$,~$\tilde D^i_t=D_t$ for $t \in [x_i,x_{i+1}]$ and $\tilde D^i_t=D_{x_{i+1}}+\phi(t-x_{i+1}) Q$ for $t \in [x_{i+1},x_{i+1}+1]$. We define $(\tilde D^0_t)_{t \in [x_0,x_1+1]}$ by $\tilde D^0_t=D_t$ for $t \in [x_0,x_1]$ and $\tilde D^0_t=D_{x_1}+\phi(t-x_1)Q$ for $t \in [x_1,x_1+1]$ and furthermore $(\tilde D^k_t)_{t \in [x_k-1,x_{k+1}]}$ by $\tilde D^k_t:=D_{x_k}+\phi(x_{k}-t)Q$ for $t \in [x_k-1,x_k]$ and $\tilde D^k_t=D_t$ for $t \in [x_k,x_{k+1}]$. 

For each path $\tilde D_t^i$ the previous part of the proof works. 
Furthermore  $$\spfl(D_t)=\sum_{i=0}^k \spfl(\tilde D^i_t) \ .$$ For notational simplicity we assume that $k=1$. 
Then by the first part of the proof
\begin{eqnarray*}
\lefteqn{\spfl(\tilde D^0_t)+\spfl(\tilde D^1_t)}\\
&=& \frac 12~\int_0^{x_1}(\frac{d}{dt}D_t)\chi'(D_t)~ dt + \frac 12~\int_{x_1}^{x_1+1} (\frac{d}{dt}\tilde D_t^0)\chi'(\tilde D^0_t)~ dt \\
&& + ~  \frac 12~\int_{x_1-1}^{x_1} (\frac{d}{dt}\tilde D_t^1)\chi'(\tilde D^1_t)~ dt + \frac 12~\int_{x_1}^1(\frac{d}{dt}D_t)\chi'(D_t)~ dt \\
&& +~ \frac 12~ \Tr(21_{\ge 0}(\tilde D^0_{x_1+1})-1-\chi(\tilde D^0_{x_1+1})) - \frac 12 \Tr(2P_0-1-\chi(D_0))\\
&& +~ \frac 12 \Tr(2P_1-1-\chi(D_1)) - \frac 12 \Tr(21_{\ge 0}(\tilde D^1_{x_1-1})-1-\chi(\tilde D^1_{x_1-1}))
\end{eqnarray*}
The assertion follows since 
$$\int_{x_1}^{x_1+1} (\frac{d}{dt}\tilde D_t^0)\chi'(\tilde D^0_t)~ dt=-\int_{x_1-1}^{x_1} (\frac{d}{dt}\tilde D_t^1)\chi'(\tilde D^1_t)~ dt$$ 
and $$\Tr(21_{\ge 0}(\tilde D^0_{x_1+1})-1-\chi(\tilde D^0_{x_1+1}))=\Tr(21_{\ge 0}(\tilde D^1_{x_1-1})-1-\chi(\tilde D^1_{x_1-1})) \ .$$

Now consider general $\chi$. We may assume that $\supp \phi=\bbbr$. Let $f \in \C_c(\bbbr)$ be an even function with $0\le f \le 1$ and $f(0)=1$ and let $f_n(x)=f(\frac xn)$ for $n \in \bbbn$. Let $\psi_n=\chi' f_n$. Let $C_n=\int_0^{\infty} \psi_n(x)~dx$ and $\chi_n(x)=\frac{1}{C_n}\int_0^x \psi_n(y)~dy$. Then $\chi_n$ is a normalizing function for each $n$. Furthermore $\chi_n$ converges to $\chi$ in $C^1(\bbbr)$ for $n \to \infty$. For $\chi_n$ the formula holds by the previous part of the proof. Since $f_n(D_t)$ is uniformly bounded and converges to the identity in $l^{\infty}(H)$ for $n \to \infty$, it follows that $(\frac{d}{dt}D_t)\psi_n(D_t)$ converges to $(\frac{d}{dt}D_t)\chi'(D_t)$ uniformly in $l^1(H)$. 

Furthermore for $x \le 0$ we have that $$1-|\chi_n(x)|=\frac{1}{C_n}\int_{-\infty}^x \psi_n(y)~dy \le \frac{1}{C_1}(1-|\chi(x)|)$$ and similarly for $x>0$. Hence for $x \in (-\infty,0)$ 
$$|1-\chi_n^2(x)| \le \frac{2}{C_1}\phi(x) \ .$$  Thus $1-\chi_n(D_i)^2$ converges to $1-\chi(D_i)^2$ in $l^1(H)$ for $i=0,1$. Furthermore $(\chi_n(D_i) +(2P_i-1))^{-1}$ converges to $(\chi(D_i) +(2P_i-1))^{-1}$ in $l^{\infty}(H)$. It follows that $$\chi_n(D_i)-(2P_i-1)=(\chi_n^2(D_i)-1)(\chi_n(D_i) +(2P_i-1))^{-1}$$ converges to $\chi(D_i)-(2P_i-1)$ in $l^1(H)$.
\end{proof}

\begin{cor} Assume that $(D_t)_{t \in [0,1]}$ is a path of selfadjoint operators with compact resolvents and common domain such that $(t \mapsto D_t) \in C^1([0,1],B(H(D_0),H))$ and $D_0$, $D_1$ are invertible and assume that there is $U \in \U(H)$ with $UD_0U^*=D_1$.

Let $\phi \in C_0(\bbbr)$ and assume that $\supp \phi$ is connected and $\phi$ is positive on the interior of its support and that $[0,1] \to l^1(H),~t \mapsto \phi(D_t)$ is continuous.

Let $\psi \in C_0(\bbbr)$ be an even non-negative function with $\supp \psi \subset (\supp \phi)^o,~\psi(0)>0$ and $(|x|+1)^p |\psi(x)| \le C\phi(x)$ where we set $p=0$ if $D_0-D_t$ is bounded for all $t \in [0,1]$, and else $p=1$. Let $$C=\int_{-\infty}^{\infty}\psi(z)~dz <\infty \ .$$ Then
$$\spfl((D_t)_{t \in [0,1]}) = \frac{1}{C} \int_0^1  \Tr (\frac{d}{dt}D_t)\psi(D_t)~ dt \ .$$
\end{cor}

\begin{proof}
As in the last part of the proof of the theorem let $f \in \C_c(\bbbr)$ be an even function with $0\le f \le 1$ and $f(0)=1$ and let $f_n(x)=f(\frac xn)$ for $n \in \bbbn$. Let $\psi_n=\psi f_n$. Let $C_n=\int_0^{\infty} \psi_n(x)~dx$ and $\chi_n(x)=\frac{1}{C_n}\int_0^x \psi_n(y)~dy$. The assertion holds since $(\frac{d}{dt}D_t)\psi_n(D_t)$ converges uniformly in $l^1(H)$ to $(\frac{d}{dt}D_t)\psi(D_t)$. 
\end{proof}

\section{Applications}

\subsection{Paths of elliptic operators} Let $M$ be a closed Riemannian manifold and $E$ a hermitian vector bundle on $M$. Let $(D_t)_{t \in [0,1]}$ be a path of elliptic symmetric differential operators of order one on $L^2(M,E)$ with invertible endpoints. Assume that the coefficients of $D_t$ depend smoothly on $t$. 
There is $p\ge 1$ such that the injection $H^2(M,E) \to L^2(M,E)$ is $p/2$-summable (with respect to any identification $H^2(M,E) \cong L^2(M,E))$. Since $$(1+D_t^2)^{-1}:L^2(M,E) \to H^2(M,E)$$ depends continuously on $t$, we have that $(1+D_t^2)^{-p/2}$ depends continuously on $t$ in $l^1(L^2(M,E))$. Let $\chi \in C^1(\bbbr)$ be a normalizing function such that there is $C>0$ with $|\chi^2(x)-1| \le C(1+x^2)^{-p/2}$ and $(|x|+1)|\chi'(x)| \le C(1+x^2)^{-p/2}$. 
Then the formula in the theorem holds for $\chi$ and $(D_t)_{t \in [0,1]}$.

Particular examples of normalizing functions fulfilling these conditions are the functions $\chi_p$ and $\chi^s$ defined in the following section.

By the remark before the proof of the theorem the formula does not depend on the choice of the Riemannian metric on $M$ and the hermitian structure on $E$, which are even allowed to vary smoothly with the parameter. Note that for the definition of $\frac{d}{dt}D_t$ on $\C(M,E)$ no choice of Riemannian metric and hermitian structure is required.

\subsection{$p$-summable operators and $\theta$-summable operators} 
We consider the spectral flow of $(D_t:=D+A_t)_{t \in [0,1]}$ where $D$ has compact resolvents, furthermore $D+A_i,~i=0,1,$ are invertible and one of the following two conditions is fulfilled:
\begin{itemize}
\item[(I)]  {\it ($p$-summable case.)} $(1+D^2)^{-1/2} \in l^{p}(H)$ for some $p\ge 1$ and $(t \mapsto A_t) \in C^1([0,1],l^{\infty}(H))$. 
\item[(II)] {\it ($\theta$-summable case.)} $e^{-s_0D^2} \in l^1(H)$ for some $s_0>0$ and $(t \mapsto A_t) 
\in  C^1([0,1],l^{\infty}(H))$ and the map $[0,1] \to l^{\infty}(H),~t \mapsto [D,A_t](D+i)^{-1}$ is continuous.
\end{itemize}

Assume (I). From the resolvent formula
$$(D +A_t \pm i)^{-1}-(D\pm i)^{-1}=-(D +A_t \pm i)^{-1}A_t(D \pm i)^{-1}$$ it follows that the maps
$[0,1] \to l^p(H),~t \mapsto (D +A_t \pm i)^{-1}$ are continuous. Hence $$[0,1] \to l^1(H),~t \mapsto (1+(D+A_t)^2)^{-p/2}$$ is well-defined and continuous. 

A particular normalizing function is $$\chi_p(x):=\frac{2}{C_p}\int_0^x (1+z^2)^{-\frac{p+1}{2}}~dz$$
with  $$C_p=\int_{-\infty}^{\infty} (1+z^2)^{-\frac{p+1}{2}}~dz \ .$$ 
There is $C>0$ such that
$$|\chi_p^2(x)-1| \le 2|(21_{\ge 0}(x)-1)-\chi_p(x)| \le C(1 +x^2)^{-p/2} \ .$$ 
Furthermore $$|\chi_p'(x)|=\frac{2}{C_p}(1+x^2)^{-\frac{p+1}{2}} \le \frac{2}{C_p}(1+x^2)^{-p/2} \ .$$

Now consider (II). Define the normalizing function $$\chi^s(x):=\sqrt{\frac{s}{ \pi}}\int_0^x e^{-sz^2}~dz \ .$$  For fixed $s>0$ we have that $|(21_{\ge 0}(x)-1)-\chi^s(x)| \le Ce^{-s x^2}$ and $|{\chi^s}'(x)|\le C e^{-sx^2}$.

In the following we show that $$[0,1] \to l^1(H),~t \mapsto e^{-s(D+A_t)^2}$$ is continuous for $s>2s_0$.

There is $C>0$ such that for small $s$ 
\begin{eqnarray*}
\|((D+A_t)^2-D^2)e^{-sD^2}\|&=&\| (DA_t+A_tD+A_t^2)e^{-sD^2}\|  \\
&=& \|([D,A_t] + 2A_t D +A_t^2)e^{-sD^2}\| \\
&\le& C s^{-\frac 12} \|R_t(D+i)^{-1}\|
\end{eqnarray*}
with $R_t:=[D,A_t] + 2A_tD +A_t^2$.
For $s>0$ the map $[0,1] \to B(H),~t \mapsto R_te^{-sD^2}$ is continuous.
Hence for $s>0$ small enough the series
$$e^{-s(D+A_t)^2}=\sum_{n=0}^{\infty} (-1)^n s^n \int_{\Delta^n}e^{-u_0sD^2}R_t e^{-u_1sD^2}R_t \dots e^{-u_nsD^2} ~du_0 \dots du_n$$
with $\Delta_n=\{(u_0, \dots, u_n) \subset (0,1)^{n+1}~|~ \sum_{i=0}^n u_i=1 \}$ converges in $B(H)$ and depends continuously on $t$.
This implies that the operator $[0,1] \to B(H), ~t \mapsto e^{-s(D+A_t)^2}$ is continuous.

Furthermore by \cite[Theorem C]{gs} for $s_1>s_0$ we have that $e^{-s_1(D+ A_t)^2} \in l^1(H)$ for any $t \in [0,1]$ and there is $\lambda(s_0,s_1) \in \bbbr$ such that
$$\|e^{-s_1(D+ A_t)^2}\|_1 \le e^{\lambda(s_0,s_1)\|A_t\|^2}\|e^{-s_0D^2}\|_1 \ .$$

Let $s>2s_1>2s_0$. For $h$ small
\begin{eqnarray*}
\lefteqn{e^{-s(D+A_t)^2}-e^{-s(D+A_{(t+h)})^2}}\\
&=&(e^{-(s-s_1)(D+A_t)^2}-e^{-(s-s_1)(D+A_{(t+h)})^2})(e^{-s_1(D+A_t)^2} + e^{-s_1(D+A_{(t+h)})^2})\\
&& +  ~e^{-(s-s_1)(D+A_{(t+h)})^2}e^{-s_1(D+A_t)^2}   - e^{-(s-s_1)(D+A_t)^2}e^{-s_1(D+A_{(t+h)})^2} \ .
\end{eqnarray*}
For $h \to 0$ the first summand on the right hand side converges to zero in $l^1(H)$ since the first factor does so in $B(H)$ and the second factor is uniformly bounded in $l^1(H)$ by the previous estimate. Since $e^{-s_1(D+A_t)^2} \in l^1(H)$ and $e^{-(s-s_1)(D+A_t)^2} \in l^1(H)$, clearly the second row on the right hand side also converges to zero in $l^1(H)$. It follows that for $s>2s_0$
$$[0,1] \to l^1(H),~ t \mapsto e^{-s(D+A_t)^2}$$ is continuous.

Summarizing, we conclude that in situation (I) the second condition of the theorem and the corollary is fulfilled for $\phi(x)=(1+x^2)^{-\frac p2}$ and in situation (II) for $\phi(x)=e^{-sx^2}$ with $s > 2s_0$. 

The theorem then implies for example the following formulas: 

In situation (I) for $C_p$ and $\chi_p$ as above
\begin{eqnarray*}
\spfl((D_t)_{t \in [0,1]})&=&(C_p)^{-1}\int_0^1  \Tr (\frac{d}{dt}D_t)(1+D_t^2)^{-\frac{p+1}{2}}~ dt\\
&& + ~\frac 12\Tr((2P_1-1)-\chi_p(D_1)) - \frac 12\Tr((2P_0-1)-\chi_p(D_0)) \ .
\end{eqnarray*} 
If $D+A_0$ and $D+A_1$ are unitarily equivalent, then by the corollary
$$\spfl((D_t)_{t \in [0,1]})=(C_{p-1})^{-1}\int_0^1  \Tr (\frac{d}{dt}D_t)(1+D_t^2)^{-p/2}~ dt \ .$$
In situation (II) for $s>2s_0$ and $\chi^s$ as above
\begin{eqnarray*}
\spfl((D_t)_{t \in [0,1]})&=&\sqrt{\frac{s}{\pi}} \int_0^1  \Tr (\frac{d}{dt}D_t)e^{-sD_t^2}~ dt\\
&& +~ \frac 12\Tr((2P_1-1)-\chi^s(D_1)) - \frac 12\Tr((2P_0-1)-\chi^s(D_0))\ .
\end{eqnarray*} 
In comparison with \cite{cp2}, where these formulas have been proven for norm continuous paths $(A_t)_{t \in [0,1]}$ , we only require strong continuity. However, we have an additional condition in the $\theta$-summable case. Also the contributions of the endpoints are different.

Another motivation for reconsidering the formulas here is that the proofs in \cite{cp1}\cite{cp2} are technically very involved. This is partially due to the fact that they include the case of Breuer--Fredholm operators affiliated to a semifinite von Neumann algebra. Another reason is that their method requires that $t \mapsto D_t(1+D_t^2)^{-\frac 12}$ is differentiable in an appropriate topology (see p.e. \cite[\S 7.2]{bcp}). This is much more difficult to prove than the differentiability of the path of resolvents in $l^p(H)$ resp. of the path $e^{-sD_t^2}$ in $l^1(H)$, which was enough in our approach.

\end{document}